\documentclass[11pt]{amsart}
\usepackage{amsmath}
\usepackage{amssymb}
\usepackage{tabularx}
\usepackage{enumerate}
\usepackage[dvipdfm]{graphicx}
\usepackage{texdraw}
\usepackage{color}

\usepackage{mathrsfs}
\usepackage{amsfonts,amssymb,amsmath}
\usepackage{epsfig}

\makeatletter
\@addtoreset{equation}{section}

\makeatother
\newtheorem{thm}{Theorem}[section]
\newtheorem{lem}[thm]{Lemma}

\newtheorem{prop}[thm]{Proposition}
\newtheorem{remark}[thm]{Remark}

\newcommand{\R}{\mathbb{R}}

\begin{document}
\title[Sharp Waves of Porous Media Equations]
{Traveling Sharp Waves of Porous Media Equations in Spatially Periodic Environment$^\S$}
\thanks{$\S$ This research was partly supported by NSFC (No. 12071299, 12471199). }

\author[B. Lou]{Bendong Lou$^{\dag}$}
\thanks{$\dag$ Mathematics and Science College, Shanghai Normal University, Shanghai 200234, China.}
\thanks{{\bf Email:} {\sf lou@shnu.edu.cn} (B. Lou)}

\begin{abstract}
We consider one dimensional porous media equations in spatially periodic environment. We will construct a periodic traveling sharp wave whose profile tends to a positive steady state at left infinity and takes zero on the right half line, with a free boundary satisfying the Darcy's law. Our method is to take the limit for a sequence of normalized solutions starting at Heaviside type of initial data. The crucial step is to give a uniform positive lower bound for the instantaneous speed of the free boundary.
\end{abstract}

\maketitle


\section{Introduction}
In this paper we consider the following one dimensional porous media equation (PME, for short) with a spatially heterogeneous reaction:
\begin{equation}\label{inhomo-PME}
u_t = (u^m)_{xx} + f(x,u)[\kappa(x)-u],\qquad x\in \R,\ t>0,
\end{equation}
where $m>1$,
\begin{equation*}\label{F1}
\left\{
\begin{array}{l}
f(x,u)\in C^2(\R^2) \mbox{ and } \kappa(x)\in C^2(\R),\ \mbox{both are $1$-periodic in }x,\ f(x,0)\equiv 0,\\ f(x,u)>0 \mbox{ for }x\in \R,\ u>\theta \mbox{ and some }\theta>0. \mbox{ Moreover}, \ \kappa(x)>\theta \mbox{ for }x\in \R.
\end{array}
\right.
\tag{F1}
\end{equation*}
Special examples of such $f$ include
\begin{enumerate}[(a)]
\item {\it monostable case}: $f(x,u)>0$ for $x\in \R$ and $u>0$;
\item {\it bistable case}: there exists $\theta>0$ such that $f(x,u)< 0$ for $u\in (0,\theta)$, and  $f(x,u)>0$ for $u>\theta$;
\item {\it combustion case}: there exists $\theta>0$ such that $f(x,u)= 0$ for $u\in [0,\theta]$, and  $f(x,u)>0$ for $u>\theta$;
\item {\it multistable case}: there exists $\theta>0$ such that $f(x,\cdot)$ changes sign finitely many  times in $(0,\theta)$, while $f(x,u)>0$ for $u>\theta$.
\end{enumerate}
The purpose of this paper is to construct traveling waves of \eqref{inhomo-PME}.
Such waves are called {\it periodic traveling waves} or {\it pulsating traveling waves} in spatially periodic reaction-diffusion equations (RDEs, for short, that is, the case where $m=1$). In the last decades, they have been extensively studied in the RDE cases (cf. \cite{BH2002, BHN2005, BH2007, DGM, Xin1991, Xin2000, Xin-book} and references therein). In heterogeneous PME cases, however, not much is known so far, even for the special monostable, bistable or combustion equations.

In spatially homogeneous cases, the PME (especially PME with multistable reactions) may have two types of traveling waves (cf. \cite{A1, GK, Vaz-book}). The profile of each of the first type is positive on the whole line, like that in RDEs. The second type, however, tends to a positive steady state at left infinity
and takes $0$ on the right side of a free boundary. It is called a {\it sharp wave} by many authors, and its existence is due to  the degeneracy of the equation. In some sense, the sharp wave is more essential, because only this type exists in monostable, bistable and combustion PMEs. In this paper, we will study the analogue for \eqref{inhomo-PME}, and call it a {\it periodic sharp wave} for simplicity. 

Note that, our equation \eqref{inhomo-PME} includes bistable and even multistable cases.
Without further conditions, rightward propagating traveling wave may not exist. For example, in the homogeneous bistable case: $f(x,u)=u(u-a)$ and $\kappa(x)\equiv 1$, rightward traveling wave exists if and only if $a\in (0,\frac12)$.
For this reason, we impose the following additional assumption:
\begin{equation*}\label{F2}
\left\{
 \begin{array}{l}
 \mbox{there exist } c_0, l_0>0,\ \varphi_0(z)\in C^2([-l_0,l_0]) \mbox{ with }\varphi_0(\pm l_0)=0,\\
  \varphi_0(z)>0 \mbox{ in }(-l_0,l_0),\  -\varphi'_0(l_0)> c_0,  \mbox{ and }\\
(\varphi_0^m)'' (z)+ c_0 \varphi'_{0}(z) + f(x,\varphi_0(z))[\kappa(x)-\varphi_0(z)]\geq 0 \mbox{ for } x\in \R, z\in [-l_0,l_0].
\end{array}
\right.
\tag{F2}
\end{equation*}
This means that $\varphi_0(x-c_0 t)$ is a rightward moving subsolution of \eqref{inhomo-PME}, and so it  guarantees the possibility of traveling waves. In Section 2, we will explain that \eqref{F2} holds automatically in monostable and combustion cases, while it holds in the bistable case under some easily verified sufficient conditions.

Our first main result is about the existence of the periodic sharp wave.
\begin{thm}\label{thm:main1}
Assume \eqref{F1} and \eqref{F2}. Then the equation \eqref{inhomo-PME} has a {\rm periodic traveling wave}
$U(x,t)\in C(\R^2)$ with
\begin{equation}\label{U=PTW}
U(x,t+T)\equiv U(x-1,t),\quad x,\ t\in \R,
\end{equation}
for some $T>0$. Moreover,
\begin{enumerate}[{\rm (i)}]
\item there exists $B(t)\in C^1(\R)$ such that
$$
U (x,t)>0 \mbox{ in } Q:=\{(x,t)\in \R^2 \mid x<B(t),\ t\in \R\}, \qquad U(x,t)\equiv 0\mbox{ in } Q^c := \R^2\backslash Q;\
$$
and $U\in C^{2+\alpha,1+\alpha/2}(Q)$ for any $\alpha\in (0,1)$;
\item there exists a positive, $1$-periodic stationary  solution $p_0(x)$ of \eqref{inhomo-PME} such that
$$
U(x,t)-p_0(x)\to 0\mbox{ as }x\to -\infty, \qquad \mbox{for each } t\in \R;
$$
\item there exists $\delta^*>0$ such that the Darcy's law holds in the following sense
$$
B'(t)= -\frac{m}{m-1} [U^{m-1}]_x (B(t)-0,t)\geq \delta^*,\qquad t\in \R.
$$
\end{enumerate}
\end{thm}

This theorem is proved by taking limit for a sequence of normalized solutions. More precisely,
we will consider the equation \eqref{inhomo-PME} with initial data
\begin{equation}\label{ini-PME}
u_0(x) = p(x)H(-x),\qquad x\in \R,
\end{equation}
where $H(x)$ is the Heaviside function and $p(x)$ is a positive stationary solution of \eqref{inhomo-PME} satisfying
\begin{equation}\label{p>varphi}
p(x)\geq \varphi_0(x),\qquad x\in [-l_0,l_0]
\end{equation}
for $\varphi_0$ in \eqref{F2}.
Denote the solution of this Cauchy problem by $u(x,t)$ and denote its unique free boundary by $b(t)$. In addition, for any positive integer $n$, define $t_n$ as the time when $u(n,t_n)=0$. Then the periodic traveling wave $U$ and its free boundary $B(t)$ in the previous theorem can be obtained as in the following theorem.

\begin{thm}\label{thm:main-limit}
Assume \eqref{F1} and \eqref{F2}. Let $u(x,t),\ b(t)$ and $t_n$ be as above.
Then
\begin{equation}\label{u-to-U}
u(x+n,t+t_n)\to U(x,t) \mbox{\ \ as\ \ }n\to \infty,
\end{equation}
in the topology $C_{loc}(\R^2)$ and $C^{2+\alpha,1+\alpha/2}_{loc}(Q)$ for any $\alpha\in (0,1)$, and
\begin{equation}\label{b-to-B}
b(t+t_n)\to B(t) \mbox{\ \ as\ \ }n\to \infty,
\end{equation}
in the topology $C_{loc}(\R)$.
\end{thm}

The above two theorems will be proved in Section 3. Our approach is the standard renormalization method, that is, taking limit for the sequence $\{u(x+n,t+t_n)\}$. We will see below that this sequence is decreasing in $n$, so, how to avoid the limit $U$ being the trivial $0$ is a big problem. For this purpose, we will give uniform positive lower bounds for the instantaneous speed $b'(t)$ and for $u(x,t)$ when $x$ lies on the left side of the free boundary and away from it. This will be the main difficulties to be solved in Section 3.

We remark that, one can also renormalize $u(x,t)$ by considering the $h$-level set (as it was done in \cite{DGM} for RDEs): for some $h>0$, set $s_n$ the time when $u(n,s_n)=h$, and consider the renormalized sequence $\{u(x+n,t+s_n)\}$. If this sequence converges to some entire solution $\tilde{U}(x,t)$, then it can be expected to be a periodic traveling wave. Of course it is non-trivial. The problem, however, is that we do not know whether $\tilde{U}(x,t)$ is a sharp wave (that is, it has a free boundary) or it is positive on the whole line $\R$. Our renormalization process can indeed give a sharp wave.

The convergence in \eqref{u-to-U} is generally not true in $L^\infty(\R)$ topology, since the equation can be a multistable one, and the sequence $\{u(x+n,t+t_n)\}$ may \lq\lq converge" to a {\it propagating terrace} (cf. \cite{DGM} for RDEs). In such cases, the sharp wave is actually the lowest component of the terrace, and so the convergence in \eqref{u-to-U} does not hold in $L^\infty(\R)$ norm.
Nevertheless, in the special cases like the monostable, bistable and combustion equations, the $L^\infty$ convergence can be true.

\begin{thm}\label{thm:L-infty-cov}
Assume \eqref{F1}, $p(x)$ in \eqref{ini-PME} is the minimal $1$-periodic stationary solution of \eqref{inhomo-PME} in the range $[\min\kappa(x), \max\kappa(x)]$. Then
\begin{equation}\label{Linfty-conv}
\|u(\cdot+n, t+t_n) - U(\cdot,t)\|_{L^\infty(\R)} \to 0 \mbox{ as }n\to \infty,
\end{equation}
locally uniformly in $t\in \R$, if one of the following holds
\begin{enumerate}[{\rm (i)}]
\item the reaction term is a monostable one as in Example {\rm (a)};
\item the reaction term is a combustion one as in Example {\rm (c)};
\item the reaction term is a bistable one as in Example {\rm (b)}, and \eqref{F2} holds with $\varphi_0(x)\leq p(x)$.
\end{enumerate}
\end{thm}

The aim of seeking for traveling waves is, roughly speaking, to use them to characterize the spreading phenomena of the solutions, see, for example, the classic references \cite{AW1, Fisher, FM, KPP} etc.
In particular, in \cite{FM} Fife and McLeod proved that, a solution of homogenous bistable RDE starting at a compactly supported initial data will converge as $t\to \infty$ to a combination of traveling waves moving opposite each other. Similar results have been obtained for homogeneous RDEs and PMEs.
In a forthcoming paper \cite{LouLu}, we will also use the periodic sharp traveling wave to characterize the
spreading solutions of the Cauchy problem of \eqref{inhomo-PME}.
%
%
%

The paper is arranged as the following. In Section 2 we give some preliminaries, including basic results one the solutions of PMEs, positive stationary solutions, sufficient conditions for \eqref{F2}, pressure equation, intersection number properties and general convergence results. In Section 3 we specify the initial value problem, give the a priori estimates, take renormalization, and prove Theorems \ref{thm:main1}, \ref{thm:main-limit}. In Section 4 we consider the $L^\infty(\R)$ convergence in Theorem \ref{thm:L-infty-cov}.

\section{Preliminaries}
In this section we give some preliminary results, including the definition of weak solutions and their properties, positive stationary solutions, discussion of the condition \eqref{F2}, pressure equation, intersection number properties and general convergence results.

\subsection{Basic Results}
Since the equation is degenerate at $u=0$, researchers generally define weak solutions as the following. For any $\tau >0$, denote $Q_\tau := \R \times (0,\tau)$. A function $u(x,t;u_0)\in C(Q_\tau)\cap L^\infty (Q_\tau)$ is called a {\it very weak solution} of
\eqref{inhomo-PME} with initial data $u_0(x)\in L^\infty(\R)$ if for any $\phi \in C^\infty_c (\overline{Q_\tau})$, there holds
\begin{equation}\label{def-very weak sol}
\int_{\R} u(x,\tau)\phi(x,\tau) dx  =  \int_{\R} u_0(x) \phi(x,0)dx +
\iint_{Q_\tau} f(x,u)[\kappa(x)-u] \phi  dx dt  + \iint_{Q_\tau} [u\phi_t + u^m \phi_{xx}] dx dt.
\end{equation}
As an extension of the definition, if $u$ satisfies \eqref{def-very weak sol} with inequality \lq\lq $\geq$" (resp. \lq\lq $\leq$"), instead of equality, for every test function $\phi\geq 0$, then $u$ is called a
{\it very weak supersolution} (resp. {\it very weak subsolution}) of the problem (cf. \cite[Chapter 5]{Vaz-book}).

By the theory of PMEs  we know the following facts.
\begin{enumerate}[(1).]
\item {\it Free boundary}. Free boundaries appear in the solution when the initial data $u_0$ is a compactly supported one. For example, if $u_0(x)=  (1-x^2)_+ := \max\{0, 1-x^2\}$, then the solution will develop two free boundaries $b_l(t)<b_r(t)$.
\item {\it Darcy's law}. Any free boundary satisfies the Darcy's law. For example, at a right free boundary $b_r(t)$, there holds
\begin{equation}\label{general-sol-Darcy}
b'_r(t) = - \frac{m}{m-1}[u^{m-1}]_x (b_r(t)-0,0),\qquad t>0.
\end{equation}
Hereinafter, $u_x(a-0,t)$ denotes the left-hand derivative of $u$ at $x=a$.
\item {\it Positivity persistence and regularity}. For any $x_1\in \R$, once $u(x_1,t_1)>0$ for some $t_1\geq 0$, then $u(x_1,t)>0$ for all $t\geq t_1$. $u(x,t)$ is classical in a neighborhood of $(x_1,t_1)$ if $u(x_1,t_1)>0$.
\end{enumerate}
These properties can be found or be derived by the standard theory of PMEs, see, for example \cite{LouZhou, Vaz-book, Wu-book}) and references therein.

\subsection{Positive Stationary Solutions}
We now consider stationary solutions of \eqref{inhomo-PME}. Denote
$$
0<\kappa_0 := \min\limits_{0\leq x\leq 1} \kappa(x) \leq \kappa^0 := \max\limits_{0\leq x\leq 1}\kappa(x).
$$

\begin{prop}\label{prop:ppss}
Assume \eqref{F1}. Then \eqref{inhomo-PME} has $1$-periodic stationary solutions $p_1(x)$ and $p_2(x)$, which are respectively the minimal and the maximal ones of such solutions in the range $[\kappa_0, \kappa^0]$.
\end{prop}

\begin{proof}
Consider the solution $u(x,t; \frac{\kappa_0 +\theta}{2})$ of \eqref{inhomo-PME} with initial data $\frac{\kappa_0 +\theta}{2}$, where $\theta>0$ is the constant in \eqref{F1}. It follows from \eqref{F1} that $u(x,t)$ exists globally and is increasing in $t$. Thus $u(x,t)\to p_1(x)$ as $t\to \infty$ pointwisely. Since $u(x,t)\geq \frac{\kappa_0 +\theta}{2}$, it is a classical solution. The convergence $u\to p_1$ actually holds in $C^{2,1}_{loc}(\R^2)$ topology. This implies that $p_1(x)$ is a positive stationary solution of \eqref{inhomo-PME}, and it is the minimal one in the range $[\kappa_0,\kappa^0]$. Moreover, since the equation is invariant when we replace $x$ by $x+1$, we see that $u(x,t)\equiv u(x+1,t)$, and so $p_1(x+1)\equiv p_1(x)$.

To obtain $p_2(x)$ we only need to take limit as $t\to \infty$ for the solution $u(x,t;\kappa^0 +1)$.
\end{proof}

Generally, one does not know whether $p_1\equiv p_2$ or not (except for the Fisher-KPP equations).
In addition, besides the stationary solutions in the range $[\kappa_0,\kappa^0]$, \eqref{inhomo-PME} may have some other stationary solutions. All of them are less than $p_2(x)$. For example, in the bistable and combustion examples (b) and (c), $u\equiv \theta$
is such a solution.

\subsection{Discussion of Condition \eqref{F2}}
We now discuss the condition \eqref{F2} and give some sufficient condition for it.
For this purpose, we first recall some results on the traveling waves of the following homogeneous PME:
\begin{equation}\label{homo-PME-eq}
u_t = (u^m)_{xx} + f_0(u),
\end{equation}
where $m>1$, $f_0\in C^2([0,\infty))$ and $f_0(0)=0$. If there is a function $q(z)\in C^2([-l,l])$ for some $l>0$ such that $q(z)>0$ in $(-l,l)$ and
\begin{equation}\label{TW-q}
(q^m)'' + cq' + f_0(q)=0,\qquad z\in [-l,l],
\end{equation}
they we call $u=q(x-ct)$ a traveling wave of the equation \eqref{homo-PME-eq}.
Note that it is only a traveling-wave-type of solution of the equation \eqref{homo-PME-eq} on the interval $x\in [ct-l,ct+l]$, but not a traveling wave of the Cauchy problem of \eqref{homo-PME-eq} since it does not necessarily satisfy the Darcy's law on the free boundaries.

By using the phase plane analysis as in \cite{A1, AW1}, one can show the following result.
\begin{prop}\label{prop:compact-TW}
Assume $f_0 \in C^2([0,\kappa_0])$ satisfies $f_0(0)=f_0(\kappa_0)=0$,  and
\begin{equation}\label{h-cond}
\int_u^{\kappa_0} r^{m-1} f_0 (r) dr >0 \mbox{ for all } u\in [0,\kappa_0).
\end{equation}
Then for any small $c>0$ and any $q_0\in (0,\kappa_0)$ close to $\kappa_0$, there exist $l (c,q_0)>0$ and $q(z;c,q_0)$ such that
\begin{equation}\label{prop-q}
q(z;c,q_0)>0 \mbox{ in } (-l(c,q_0),l(c,q_0)),\quad q(z;c,q_0)=0 \mbox{ for }z=\pm l(c,q_0),
\end{equation}
$q$ solves \eqref{TW-q} in $[-l(c,q_0),l(c,q_0)]$, and
\begin{equation}\label{homo-sub-TW}
c< -\frac{m}{m-1} [q^{m-1}]_z (l(c,q_0)-0;c,q_0)
\end{equation}
when $c>0$ is small.
\end{prop}

\noindent
By \eqref{homo-sub-TW}, $q(x-ct;c,q_0)$ is indeed a traveling-wave-type of subsolution of the Cauchy problem of \eqref{homo-PME-eq}. Assume
\begin{equation*}\label{F3}
\left\{
 \begin{array}{l}
 f(x,u)[\kappa(x)-u]\geq f_0(u) \mbox{ for } x \mbox{ and } u\in [0,\kappa_0], \\
 \mbox{where } f_0 \mbox{ satisfies the conditions in Proposition } \ref{prop:compact-TW}.
\end{array}
\right.
\tag{F3}
\end{equation*}
Under this assumption, it is easily seen that \eqref{F2} holds for $\varphi_0(z)=q(z;c,q_0)$, that is, \eqref{F3} is a sufficient condition for \eqref{F2}. On the other hand, \eqref{F3} is easy to be verified in many cases such as Examples (a)-(c) in Section 1. Precisely, Examples (a) and (c) satisfy \eqref{F3} automatically: we only need to choose
$$
f_0(u)=0\mbox{ for }u\in [0,\theta_1],
\qquad f_0(u) = \delta_1 (u-\theta_1)(\kappa_0 -u) \mbox{ for }u\in [\theta_1, \kappa_0],
$$
where $0<\delta_1 \ll 1$,  $\theta_1 := \frac12 \kappa_0$ in Example (a), and choose
$$
f_0(u)=0\mbox{ for }u\in [0,\theta_2], \qquad f_0(u) = \delta_2 (u-\theta_2)(\kappa_0 -u) \mbox{ for }u\in [\theta_2, \kappa_0],
$$
where $0<\delta_2 \ll 1$,  $\theta_2 := \frac{\theta +\kappa_0}{2}$ in Example (c). Example (b) also satisfies \eqref{F3} if
$$
\int_u^\theta \Big[ \min\limits_{x\in [0,L]} f(x,r) \Big](\kappa^0-r)r^{m-1} dr + \int_{\theta}^{\kappa_0} \Big[ \min\limits_{x\in [0,L]} f(x,r) \Big](\kappa_0-r)r^{m-1} dr >0,\quad u\in [0,\theta).
$$

\subsection{Intersection Number}
The so-called zero number argument (cf. \cite{Ang}) developed in 1980s is a powerful tool in the study of asymptotic behavior of solutions of one-dimensional parabolic equations.
In this paper we will use the PME version that was developed in \cite{LouZhou}, not only for the asymptotic study but also for the a priori estimates.
To state the argument more precisely, it is convenient to consider the pressure function
\begin{equation}\label{def-pressure}
v(x,t):= \frac{m}{m-1} u^{m-1}(x,t)
\end{equation}
instead of the density function $u$.  Hereinafter, we will use $v,\ v^{(i)},\ v_i,\ V,\ \cdots$ to
denote the corresponding pressure function of $u,\ u^{(i)},\ u_i,\ U,\ \cdots$.
If $u$ solves \eqref{inhomo-PME}, then $v$ solves
\begin{equation}\label{v-PME}
v_t = (m-1)vv_{xx} + v_x^2 + g(x,v), \qquad  x\in \R,\ t>0,
\end{equation}
where
\begin{equation}\label{def-G}
g(x,v):= m \Big( \frac{(m-1)v}{m}\Big)^{\frac{m-2}{m-1}} f
\left(x, \Big( \frac{(m-1)v}{m}\Big)^{\frac{1}{m-1}} \right)\Big[ \kappa(x) -\Big( \frac{(m-1)v}{m}\Big)^{\frac{1}{m-1}}\Big].
\end{equation}
We say that $v$ is a positive solution of the Cauchy problem of \eqref{v-PME} with free boundaries $l(t)<r(t)$ in the time interval $(t_1,t_2)$, we mean that not only $v$ is positive and solves \eqref{v-PME} in $(l(t),r(t))$, but also the free boundaries satisfy the Darcy's law
\begin{equation}\label{v-Darcy}
l'(t) = -v_x (l(t)+0,0),\qquad r'(t)=-v_x(r(t)-0,t),\qquad t\in (t_1,t_2).
\end{equation}

First, we specify some special relationships between two solutions. For $i=1,2$, assume $l^{(i)}(t) <r^{(i)}(t)$ are continuous functions for $t\in [t_1, t_2)$,
$$
Q^{(i)}:= \{(x,t)\mid l^{(i)}(t)<x<r^{(i)}(t),\ t_1<t<t_2\},
$$
$v^{(i)}(x,t)\in C(\R\times [t_1, t_2)) \cap C^{2,1}(Q^{(i)})$ is a solution of \eqref{v-PME}. Denote
$$
l(t) := \max\{l^{(1)}(t), l^{(2)}(t)\},\qquad r(t):=\min\{r^{(1)}(t), r^{(2)}(t)\},\qquad t\in [t_1,t_2).
$$
When $l(t)< r(t)$ we use the following notations.

\medskip
\noindent
\underline{\it Case 1}. 
Denote $v^{(1)}(\cdot,t_0) \vartriangleright v^{(2)}(\cdot,t_0)$, if there exists $x_0\in (l(t_0), r(t_0))$ such that
$$
v^{(1)}(x,t_0)> v^{(2)}(x,t_0) \mbox{ for } l(t_0)\leq x<x_0,\qquad v^{(1)}(x,t_0)<v^{(2)}(x,t_0) \mbox{ for } x_0 < x \leq r(t_0).
$$
(Note that this notation is a little different from the concept {\it $v^{(1)}$ is steeper than $v^{(2)}$} in \cite{DGM}. The later will be used as it was in Subsection 3.3.)

\medskip
\noindent
\underline{\it Case 2}. 
Denote $v^{(1)}(\cdot,t_0) \succapprox v^{(2)}(\cdot,t_0)$, if $l^{(1)}(t_0)\leq l^{(2)}(t_0) < r^{(2)}(t_0) =r^{(1)}(t_0)$, or, $l^{(1)}(t_0)=l^{(2)}(t_0) < r^{(2)}(t_0) \leq r^{(1)}(t_0)$, and
$$
v^{(1)}(x,t_0)> v^{(2)}(x,t_0) \mbox{ for } l(t) < x <  r(t).
$$

\medskip
\noindent
\underline{\it Case 3}. 
Denote $v^{(1)}(\cdot,t_0) \succ v^{(2)}(\cdot,t_0)$, if $l^{(1)}(t_0)<l^{(2)}(t_0) < r^{(2)}(t_0) <r^{(1)}(t_0)$, and
$$
v^{(1)}(x,t_0) > v^{(2)}(x,t_0) \mbox{ for } l(t) < x < r(t).
$$

\medskip
In the next section 
we actually consider such relationships among three solutions. More precisely, let $v_i(x,t)\in C(\R\times [t_1, t_2)) \cap C^{2,1}(Q_i)$\ $(i=1,2,3)$ be three solutions of the Cauchy problem of \eqref{v-PME}, where
$$
Q_i:= \{(x,t)\mid l_i(t)<x<r_i(t),\ t_1<t<t_2\},
$$
and $l_i(t)<r_i(t)$ are free boundaries of $v_i$.
Assume also that $l_1(t) \equiv l_2(t)\equiv -\infty$, other free boundaries are bounded functions.
%
%
%
The original {\it zero number argument} (cf. \cite{Ang}) says, roughly, that the zero number of a solution of a {\it one-dimensional linear uniformly} parabolic equation is decreasing. Correspondingly, the {\it intersection number} of two solutions of a nonlinear uniformly parabolic equation is also decreasing (since the difference function of these two solutions solves a linear equation). In \cite{LouZhou} the authors extended this property from uniformly parabolic equations to degenerate ones, like the PMEs. As a consequence we have the following result (see details in \cite{LouZhou}).

\begin{prop}\label{prop:inter}
\begin{enumerate}[{\rm (i)}]
\item If $v_1(\cdot ,t_1)\vartriangleright v_2(\cdot,t_1)$, then there exist $s_1, s_2$ with $t_1<s_1\leq s_2 \leq t_2$ such that $v_1(\cdot,t)\vartriangleright v_2(\cdot,t)$ for $t\in [t_1, s_1)$, $v_1(\cdot,t) \succapprox v_2(\cdot,t)$ for $t\in [s_1,s_2]$, and $v_1(\cdot,t)\succ v_2(\cdot,t)$ for $t\in (s_2, t_2)$.

\item If $r_1(t_1) = l_3(t_1)$, then either $r_1(t)=l_3(t)$ for all $t\in [t_1,t_2)$, or, there exist $s_3, s_4, s_5$ with $t_1\leq s_3 < s_4 \leq s_5 \leq t_2$ such that $r_1(t)=l_3(t)$ for $t\in [t_1,s_3]$, $v_1(\cdot,t)\vartriangleright v_3(\cdot,t)$ for $t\in (s_3, s_4)$, $v_1(\cdot,t) \succapprox v_3(\cdot,t)$ for $t\in [s_4,s_5]$, and $v_1(\cdot,t)\succ v_3(\cdot,t)$ for $t\in (s_5, t_2)$.
%
%
\end{enumerate}
\end{prop}

\noindent
In (i), the relationship $v_1(\cdot,t) \succapprox v_2(\cdot,t)$ for $t\in [s_1,s_2]$ might be true even if $s_1 <s_2$, since we have now Hopf boundary lemma for PMEs. In (ii), $r_1(t)=l_3(t)$ might be true in a period after $t_1$, since the waiting times of $r_1(t)$ and $l_3(t)$ at $r_1 (t_1)$ may be positive.
%
%
%


\section{Construction of Periodic Sharp Wave}
In this section we construct the periodic sharp traveling wave of the equation \eqref{inhomo-PME}, and prove Theorems \ref{thm:main1} and \ref{thm:main-limit}.
The idea of the proof is using the renormalization method as in \cite{DGM} to construct an entire solution, and then show that this entire solution in nothing but the desired traveling wave.
From now on, we always assume \eqref{F1} and \eqref{F2}.

\subsection{Cauchy Problems}
Let $p_2(x)$ be the positive stationary solution in Proposition \ref{prop:ppss}. Since $\varphi_0(x-c_0 t)$ in \eqref{F2} is a very weak subsolution of \eqref{inhomo-PME}, by comparison we have
\begin{equation}\label{p2>varphi}
p_2(x) > \varphi_0 (x-c_0 t-x_0),\qquad x,\ x_0,\ t\in \R.
\end{equation}
When we use the pressure functions, $p_2(x)$ corresponds to the maximal $1$-periodic stationary solution
(denoted by $q_2(x)$) of \eqref{v-PME}, and $\varphi_0(x-c_0 t)$ corresponds to a traveling-wave-type subsolution $\psi (x-c_0 t)$ of \eqref{v-PME}, where
$$
q_2 (x) := \frac{m}{m-1}p_2^{m-1}(x) \mbox{ for }  x\in \R,\qquad
\psi (z) := \frac{m}{m-1}\varphi_0^{m-1}(z) \mbox{ for }z\in [-l_0,l_0].
$$
Denote by $H(x)$ the Heaviside function:
$$
H(x)=1\mbox{ for }x\geq 0,\qquad H(x)=0\mbox{ for } x< 0.
$$
For each integer $k$, we consider the Cauchy problem
\begin{equation}\label{uk-p-v}
\left\{
\begin{array}{ll}
v_t = (m-1)vv_{xx} + v_x^2 + g(x,v), & x\in \R,\ t>0,\\
v(x,0)= v_{0k}(x):= H(k-x)q_2(x), & x\in \R.
\end{array}
\right.
\end{equation}
It follows from \cite{LouZhou, Sacks, Vaz-book, Wu-book} etc. that this problem has a unique global solution, denoted by $v(x,t;k)$, which has a unique right free boundary, denoted by $b(t;k)$, such that
$$
v(x,t;k)>0 \mbox{ in } Q_k := \{(x,t)\mid x<b(t;k),\ t> 0\},\qquad
v(x,t;k)=0 \mbox{ in } (\R\times (0,\infty))\backslash Q_k,
$$
$v\in C(\R\times (0,\infty))\cap C^{2+\alpha,1+\alpha/2}(Q_k)$ for any $\alpha \in (0,1)$, and the Darcy's law holds:
\begin{equation}\label{bn-Darcy}
b'(t;k)= - v_x (b(t;k)-0,t;k)\geq 0,\qquad t>0.
\end{equation}
(In the next subsection we will show that $b'(t;k)$ has a positive lower bound. At present, however, we only know by the standard theory of PMEs that $b(t;k)$ is strictly increasing in $t$, see, for example \cite[Theorem 2.7]{LouZhou}.) Clearly, for any integers $k$ and $j$, there holds
\begin{equation}\label{k=j}
v(x,t;k)\equiv v(x-k+j,t;j),\qquad b(t;k)\equiv b(t;j) + k-j.
\end{equation}
By \eqref{F2}, $\psi (x-c_0 t -k +l_0)$ is a subsolution of \eqref{uk-p-v} and so
\begin{equation}\label{varphi<vn}
\left\{
\begin{array}{l}
\psi(x-c_0 t -k +l_0) < v(x,t;k) \mbox{ for }x\in [c_0 t + k -2 l_0, c_0 t +k ],\ t>0,\\
c_0 t +k  <b(t;k) \mbox{ for }t>0.
\end{array}
\right.
\end{equation}
This implies that $b(t;k)\to \infty$ as $t\to \infty$.

\subsection{A Priori Estimates}\label{subsec:a-priori}
For any positive integer $n$, denote by $t_n$ the time when $v(n,t_n;0)=0$, or, equivalently, $b(t_n;0)=n$, and set
\begin{equation}\label{def:v-n}
v_n(x,t):= v(x+n,t+t_n;0),\quad b_n(t):= b(t+t_n;0)-n,\qquad  x\in \R,\ t>-t_n.
\end{equation}
Then $v_n(0,0)=0$, and
\begin{equation}\label{Darcy-bn}
b'_n(t) = -v_{nx}(b_n(t)-0,t),\qquad t>-t_n.
\end{equation}
Our purpose is to show the following convergence results:
\begin{equation}\label{conv-bn-B}
b_n(t)\to B(t) \mbox{ as }n\to \infty, \qquad \mbox{ in the topology of } C_{loc}(\R),
\end{equation}
and
\begin{equation}\label{conv-vn-V}
v_n(x,t)\to V(x,t) \mbox{ as }n\to \infty,\qquad \mbox{in the topology of } C_{loc}(\R^2) \cap C^{2+\alpha,1+\alpha/2}_{loc} (Q),
\end{equation}
where $Q:=\{(x,t)\mid x<B(t),\ t\in \R\}$ is the domain in which $V$ is positive.
Then $V(x,t)$ with free boundary $B(t)$ will be the desired periodic sharp traveling wave.
In order to show the convergence in \eqref{conv-bn-B} and \eqref{conv-vn-V}, we need some uniform-in-time a priori estimates. Due to the degeneracy of the equation, these estimates are complicated. We will specify them in several steps.

%
%
%
%

\medskip
\noindent
\underline{\it Step 1. Upper bound of $b'(t;0)$}. For any given $y>0$, define $t_y$ by $b(t_y;0)=y$. Then
$$
b(t_y;0)= y = b(t_{y+1};0) - 1 = b(t_{y+1};-1).
$$
Set $s_y := t_{y+1}-t_y$, we claim that $b(s_y; -1)>0$. In fact, if $b(s_y; -1)\leq 0$, then $v(x,s_y;-1)< v(x,0;0)$, and so by comparison we have $v(x,t;0)\succ v(x,t+s_y; -1)\equiv v(x+1,t+s_y; 0)$ for all $t>0$.
This leads to a contradiction at $t=t_y$. Using this claim and using Proposition  \ref{prop:inter} we can conclude that
\begin{equation}\label{v-sharper-vt+s}
v(x,t;0) \vartriangleright v(x,t+s_y;-1) \mbox{\ \ for\ } 0<t< t_y,
\end{equation}
and
\begin{equation}\label{exact-greater}
v(x,t_y;0) \succapprox v(x,t_{y+1};-1)\equiv v(x+1,t_{y+1};0).
\end{equation}
So,
\begin{equation}\label{speed-comp}
b'(t_y;0)=-v_x (y-0,t_y;0) \geq -v_x(y+1-0,t_{y+1};0)=b'(t_{y+1};0).
\end{equation}
(The strict inequality may not hold since we have no Hopf boundary lemma for PMEs.)
This means that $b$ moves slower at $y+1$ than it does at $y$. Hence,
\begin{equation}\label{def-sn}
s_n := t_{n+1}-t_n \mbox{ is increasing in }n,
\end{equation}
and $t_n\to \infty$ as $n\to \infty$. For any $y>0$, denote its integer part by $[y]$ and its fractional part by $\langle y \rangle$. When $[y]\geq 2$, by \eqref{speed-comp} we have
\begin{eqnarray*}
0 & \leq &  b'(t_y;0) \leq b'(t_{y-1};0)\leq \cdots \leq b'(t_{\langle y \rangle +1};0)  \\
& \leq & \bar{c} :=  \max\limits_{z\in [1,2]}  b'(t_z ;0) = \max\limits_{z\in [1,2]} \big[ -v_x(z-0,t_{z};0)\big].
\end{eqnarray*}
For any given $T_1 >0$, there holds $t_n - T_1 > t_2$ when $n$ is sufficiently large. Thus, for each large $n$ and each $t\in [-T_1,T_1]$, there exists a unique $y>2$ such that $t_y = t+t_n > -T_1 + t_n >t_2$, and so
$$
b'_n(t) = b'(t+t_n ;0) = b'(t_y; 0)\in [0, \bar{c}].
$$
By the Ascoli-Arzela lemma, a subsequence $\{b_{n_i}(t)\}$ of $\{b_n(t)\}$ converges as $i\to \infty$ to
a continuous function $B(t)$ in $C([-T_1,T_1])$. Using the Cantor's diagonal argument, we can find a subsequence of $\{b_{n_i}\}$ and a continuous function in $\R$, denote them again by $\{b_{n_i}\}$ and $B(t)$, such that, $B(0)=0$ and
\begin{equation}\label{Bni-to-B}
b_{n_i}(t) \to B(t) \mbox{ as }i\to \infty,\qquad \mbox{ in the topology of } C_{loc}(\R).
\end{equation}

\medskip
\noindent
\underline{\it Step 2. Lower bound of $v_x$ in the left neighborhood of the free boundary}.
For any given $y>3$, write $j :=[y-2]$. When $s> t_y - t_{y-j}$, as proving the claim in Step 1, one can show that $b(s;0)>j$.  Using this fact, and using a similar argument as proving \eqref{v-sharper-vt+s} and \eqref{exact-greater} we can show that
$$
v(x-j,t;0) \equiv v(x,t;j) \vartriangleright  v(x,t_y;0)\quad \mbox{\ \ for any } 0< t < t_{y-j},
$$
and
\begin{equation}\label{steeper-11}
v(x-j,t_{y-j};0) \equiv v(x,t_{y-j};j) \succapprox v(x,t_y;0).
\end{equation}
This implies that, for any $x'\in [y-1,y)$, there exists $t'\in [t_{y-j-1}, t_{y-j})$ such that the graph of $v(x,t';j)\equiv v(x-j,t';0)$ contacts that of $v(x,t_y;0)$ at exactly one point $(x',v(x',t_y;0))$, and
$$
v_x(x',t';j) \equiv v_x(x'-j,t';0) \leq v_x (x',t_y;0).
$$
Combining together with $v_x(y-j-0,t_{y-j};0) \leq v_x (y-0,t_y;0)$ (by \eqref{steeper-11}) we have
\begin{equation}\label{gradient-upper}
\min\limits_{x'\in [y-1,y]} v_x (x',t_y;0) \geq \min\limits_{{\tiny
\begin{array}{c}
x'\in [y-1,y]\\
t'\in [t_{y-j-1},t_{y-j}]
\end{array}
}
}
v_x(x'-j,t';0)  \geq -C_1:=  \min\limits_{
{\tiny
\begin{array}{c}
z\in [1,3]\\
t\in [t_1,t_3]
\end{array}
}}
v_x(z,t;0).
\end{equation}
This give the uniform-in-time lower bound for $v_x$ in the interval $[b(t;0)-1, b(t;0)]$.

As a consequence of \eqref{gradient-upper}, we have the following estimate for $v$:
\begin{equation}\label{est-near-bt}
0\leq v(x,t;0)= -\int_{x}^{b(t;0)} v_x(x,t;0)dx \leq C_1 [b(t;0)-x],\qquad b(t;0)-1 \leq x\leq b(t;0),\ t>t_3.
\end{equation}
Note that the above estimates hold only in $[b(t;0)-1, b(t;0)]$. We yet do not know the negative upper bound of $v_x$ in this interval and the positive lower bound of $v$ when $x$ is far away from $b(t;0)$, which are important when we show the limit of $v_n$ is a non-trivial one.

\medskip
\noindent
\underline{\it Step 3. Positive lower bound of the average speed}. From \eqref{speed-comp} we know that, for any positive integer $n$, it takes more time for $b(t;0)$ to cross the interval $[n+1,n+2]$ than to cross the interval $[n,n+1]$:
\begin{equation}\label{time-comp-n}
s_n :=  t_{n+1} -t_n  \leq s_{n+1} := t_{n+2} -t_{n+1}.
\end{equation}
Therefore, the average speed $\bar{c}_n:= \frac{1}{s_n}$ of $b(t;0)$ in the interval $[n,n+1]$ is decreasing in $n$. We now show that $\bar{c}_n$ has a positive limit:
\begin{equation}\label{barcn-lower-bound}
\bar{c}_n \searrow \frac{1}{T},\quad \mbox{or, equivalently}\quad s_n \nearrow T \mbox{\ \ as } n\to \infty,
\end{equation}
for some $T\in (0,\frac{1}{c_0})$, where $c_0$ is the speed of the traveling wave $\varphi(x-c_0 t +l_0)$ in \eqref{F2}. In fact, if $\bar{c}_n <c_0-\delta$ for some small $\delta>0$ and all large $n$ (to say, for $n\geq N$). Then by \eqref{varphi<vn} we have
$$
c_0 t_n  <b(t_n;0)=n-N + b(t_N;0) <(c_0 -\delta) (t_n-t_N) + b(t_N;0).
$$
This is a contradiction when $n$ is sufficiently large. This proves \eqref{barcn-lower-bound}.

\medskip
\noindent
\underline{\it Step 4. Positive lower bound of instantaneous speed $b'(t;0)$}. We will show the following lemma, which is crucial in our approach.
\begin{lem}\label{lem:positive-b'}
There exists a function $\delta_0(x)\in C^2((-\infty, 0])$ with
\begin{equation}\label{def-delta-0}
\delta_0(x) >0 \mbox{ in } (-\infty, 0), \qquad  \delta_{0*}:= \liminf\limits_{x\to -\infty} \delta_0(x)>0,\qquad \delta_0(0)=0,\qquad \delta^* := - \delta'_0(0-0)>0
\end{equation}
such that
\begin{equation}\label{v>delta0}
\left\{
\begin{array}{l}
v(x,t;0) \geq \delta_0 (x-b(t;0)),\qquad x\leq b(t;0),\ t>t_2,\\
b'(t;0) = - v_x (b(t;0)-0,t;0) \geq \delta^* = -\delta'_0(0) ,\qquad t>t_2.
\end{array}
\right.
\end{equation}
\end{lem}

\begin{proof}
Since both $f(x,u)$ and $\kappa(x)$ are of $C^2$, there exists $K>0$ such that $f(x,u)[\kappa(x) -u]\leq Ku$ for $u\geq 0$. Choose $\delta_*>0$ small such that
\begin{equation}\label{choice-delta1}
2 \delta_*^{1/2} \Big(3T +1\Big)^{1/2} e^{\frac{(m-1)K(3T+1)}{2}} <1,
\end{equation}
where $T$ is the upper bound of $s_n$ in \eqref{barcn-lower-bound}.
For any $x_0\in \R$, we consider the function
\begin{equation}\label{def-bar-u}
\bar{u} (x,t) := A e^{K(t+1)} \left( \delta_* -\frac{(x-x_0)^2}{(t+1)e^{(m-1)K(t+1)}}\right)^{\frac{1}{m-1}}_+, \qquad x\in \R,\ t>0,
\end{equation}
with $ A :=\left(\frac{m-1}{4m}\right)^{\frac{1}{m-1}}$. A direct calculation shows that
$$
\bar{u}_t - (\bar{u}^m)_{xx} - f(x,\bar{u})[\kappa(x) -\bar{u}] \geq
\bar{u}_t - (\bar{u}^m)_{xx} - K \bar{u} \geq 0,\qquad |x-x_0|\leq \rho(t),\ t>0,
$$
where
$$
\rho(t):= \delta_*^{1/2}(t+1)^{1/2} e^{\frac{(m-1)K(t+1)}{2}},\qquad t>0.
$$
So $\bar{u}$ is a supersolution, which means that
$$
\bar{u}(x,t) \geq \hat{u}(x,t;x_0),\qquad x\in [\hat{l}(t),\hat{r}(t)],\ t>0,
$$
for the solution $\hat{u}$ of the following problem
$$
\left\{
\begin{array}{ll}
\hat{u}_t = (\hat{u}^m)_{xx} + f(x,\hat{u})[\kappa(x)-\hat{u}], & x\in \R, t>0,\\
\hat{u}(x,0) = \bar{u}(x,0), & x\in \R,
\end{array}
\right.
$$
and $\hat{l}(t)=\hat{l}(t;x_0),\ \hat{r}(t)=\hat{r}(t;x_0)$ the left and right free boundaries of $\hat{u}(\cdot,t;x_0)$, respectively.
Then the right free boundary $\bar{r}(t)$ of $\bar{u}$ and that of $\hat{u}$ satisfy
$$
\hat{r}(t)\leq \bar{r}(t) = x_0 +\rho (t),\qquad t>0.
$$

Denote by $\hat{v}(x,t;x_0):= \frac{m}{m-1}\hat{u}^{m-1}(x,t;x_0)$ the pressure function corresponding to $\hat{u}$. For any $n\geq 2$, $y\in [n,n+1)$ and some suitably chosen $x_0$, we will compare $v(x,t+t_{n-2};0)$ and $\hat{v}(x,t;x_0)$ on the left side of $y$.

First, we choose $x_0=y-\rho(0)$, then the support of $\hat{v}(x,0)$ is $[y-2\rho(0), y]$. Since the waiting time for the right free boundary $\hat{r}(t)$ of $\hat{v}(x,t;x_0)$ is zero, $\hat{r}(t)$ crosses the point $y$ immediately at $t=0$. Of course, it is earlier than the free boundary $b(t+t_{n-2};0)$ of $v(x,t+t_{n-2};0)$.

Next, we choose $x_0\in [n-2+\rho(0),n-1+\rho(0)]$. Then by the choice of $\delta_*$ we have
$$
\hat{r}(3T;x_0) \leq x_0 + \rho(3T) \leq  n-1 +\rho(0) + \rho(3T) < n\leq y.
$$
This means that $\hat{r}(t;x_0)$ will use time more than $3T$ to cross the point $y$. On the other hand $b(t+t_{n-2};0)$ will use time $t_y - t_{n-2}\leq s_{n-2}+s_{n-1}+s_n\leq 3T$ to cross $y$,
earlier than $\hat{r}(t;x_0)$. Since the support of $v(x,t_{n-2};0)$ lies on the left of that of $\hat{v}(x,0;x_0)$, by Proposition \ref{prop:inter} we see that
\begin{equation}\label{b-faster-r1}
v(x,t_y;0) \succ \hat{v}(x,t_y - t_{n-2};x_0).
\end{equation}

Consequently, there exists $x_0^y\in (n-1+\rho(0), y-\rho(0))$ such that \eqref{b-faster-r1} holds for all $x_0\in [n-2+\rho(0), x_0^y)$, while when we choose $x_0=x_0^y$, there holds
\begin{equation}\label{v>v1-bdry0}
v(x,t_y;0)\succapprox \hat{v}(x,t_y-t_{n-2};x_0^y),\qquad
v(y,t_y;0)= \hat{v}(y,t_y-t_{n-2};x_0^y)=0.
\end{equation}

Now we give lower estimates for $v(x,t_y;0)$ in $I_1(y)\cup I_2(y) \cup I_3(y)$, where
$$
I_1(y):= [y-\rho(0),y],\qquad I_2(y):= [x_0^y, y-\rho(0)],\qquad I_3(y):= [x_0^y-1, x_0^y).
$$
First, we define
$$
\tilde{\delta}_1 (x) := \frac12 \min\limits_{0\leq x_0\leq 1} \min\limits_{0\leq t\leq 3T} \hat{v}(x + \hat{r}(t;x_0),t;x_0),\qquad -\rho(0)\leq x \leq 0,
$$
and take $\delta_1(x)$ as a smoothen function of $\tilde{\delta}_0(x)$ such that
\begin{equation}\label{delta-0-property}
\delta_1(x) \leq 2 \tilde{\delta}_1 (x),\qquad \delta_1(x)>0 \mbox{ in } [-\rho(0),0),\qquad \delta'_1(0)<0.
\end{equation}
Then \eqref{v>v1-bdry0} implies that
\begin{equation}\label{v>v1-bdry}
v(x,t_y;0)\succapprox \delta_1 (x-y),\qquad x\in I_1(y).
\end{equation}
Next, we consider the case where $x\in I_2(y)$. In this case, by \eqref{v>v1-bdry0} we have
\begin{equation}\label{v>v1-near-bdry}
v(x,t_y;0) \geq  \min\limits_{z\in I_2(y)} \hat{v}(z,t_y-t_{n-2};x_0^y) \geq
\delta_2:= \inf\limits_{y>2} \min\limits_{z\in I_2(y)} \hat{v}(z,t_y -t_{n-2}; x_0^y).
\end{equation}
Since $t_y -t_{n-2}$ is a finite time in $(0, 3T)$ we see that $\delta_2$ is a positive constant independent of $y$. (One can prove this point by a contradiction argument.)
Finally, we consider the case where $x\in I_3(y)$. Denote
$$
D_3 := \{(z,t)\mid x_0^y-1 \leq z\leq x_0^y , 0\leq t \leq 3 T\},\qquad
D'_3 := \{(z,t)\mid 0 \leq z\leq 1, 0\leq t \leq 3 T\}.
$$
For any $x_0\in [x_0^y -1, x_0^y)$, by \eqref{b-faster-r1} we have
\begin{equation}\label{v>v1-middle}
v(x_0,t_y;0)  \geq  \min\limits_{0\leq t\leq 3T} \hat{v}(x_0,t;x_0) \geq \min\limits_{(z,t)\in D_3 } \hat{v}(z,t;z) = \delta_3 := \min\limits_{(z,t)\in D'_3 } \hat{v}(z,t;z) >0.
\end{equation}
Combining \eqref{v>v1-bdry}, \eqref{v>v1-near-bdry} together with \eqref{v>v1-middle} we conclude that, for any $y>2$, there holds
\begin{equation}\label{v>tilde-delta}
v(x,t_y;0) \succapprox \delta_0(x-y),\qquad x\in [y-\rho(0)-1,y],
\end{equation}
if we define $\delta_0(x)$ by
$$
\delta_0(x) := \left\{
\begin{array}{ll}
\delta_1(x), & x\in [-\rho(0), 0],\\
\mbox{increasing } C^2 \mbox{ function}, & x\in [-\rho(0)-\frac12, -\rho(0)],\\
\delta_4:= \frac12 \min\{\delta_1(-\rho(0)), \delta_2, \delta_3\}, & x\in [-\rho(0)-1, -\rho(0)-\frac12].
\end{array}
\right.
$$

Finally, by the definition of $v(x,t;k)$ and by comparison, for any positive integer $k$, we have
$$
v(x,t;0)\geq v(x,t;-k)\equiv v(x+k,t;0),\qquad x\in \R,\ t>0.
$$
Thus, the lemma is proved if we extend $\delta_0(x)$ in $(-\infty, y-\rho(0)-1]$ as $\delta_4$.
\end{proof}

\subsection{Renormalization, Proof of Theorems \ref{thm:main1} and \ref{thm:main-limit}}
Based on the uniform-in-time a priori estimates in the previous subsection, we now can consider the convergence in the renormalization sequence $\{v_n\}$ and $\{b_n\}$ defined by \eqref{def:v-n}.

\medskip
\noindent
\underline{\it Step 1. Convergent subsequence}. For large $n$, from the results in the previous subsection we know the following facts:
\begin{enumerate}[(a).]
\item $v_n(x,t)\geq \delta_0 (x-b_n(t))$ for $x\leq b_n(t),\ t>-t_n$;
\item $v_n(x,t)\geq \delta_4$ for all $x\leq b_n(t)-\rho(0),\ t>-t_n$;
\item $0\leq v_n(x,t)\leq C_1 [b_n(t)-x]$ for $x\in [b_n(t)-1, b_n(t)],\ t>-t_n$.
\end{enumerate}
It follows from (a), (b) and the standard parabolic theory that, for any compact domain $K\subset Q:= \{(x,t)\mid x<B(t),\ t\in \R\}$ (where $B(t)\in C(\R)$ is the limit function in \eqref{Bni-to-B}), and any $\alpha'\in (0,1)$, there holds
$$
\|v_{n_i}(x,t)\|_{C^{2+\alpha', 1+\alpha'/2}(K)}\leq C(\alpha',K),
$$
for the same subsequence $\{n_i\}$ as in \eqref{Bni-to-B}. Hence, the sequence $\{v_{n_i}\}$ has a convergent subsequence which converges in the topology $C^{2+\alpha,1+\alpha/2}(K)$ ($\alpha\in (0,\alpha')$) to some function $V(x,t)$.
Using Cantor's diagonal argument, for any $\alpha\in (0,1)$, we can find a function $V(x,t)\in C^{2+\alpha,1+\alpha/2}(Q)$ and a subsequence of $\{n\}$ (denote it again by $\{n_i\}$) such that
\begin{equation}\label{Vni-to-V-C2}
v_{n_i}(x,t)\to V(x,t) \mbox{ as }i\to \infty,\qquad \mbox{in the topology of } C^{2+\alpha,1+\alpha/2}_{loc}(Q).
\end{equation}
In addition, $V(0,0)=0$ since $v_n(0,0)=0$.

In what follows, we extend $V$ to be zero in the domain $\R^2 \backslash Q$, and write the extended function (which is defined over $\R^2$) as $V$ again.

\medskip
\noindent
\underline{\it Step 2. To show that $V$ is a very weak entire solution and $B$ is its free boundary}.
Taking $n=n_i$ in (a) and (c), and taking limits as $i\to \infty$ we have
\begin{equation}\label{V-is-large-bdry}
\delta_0(x-B(t))\leq V(x,t)\leq C_1[B(t)-x] \mbox{ for } x\in [B(t)- 1, B(t)],\ t\in \R.
\end{equation}
This implies that $V(x,t)\in C(\R^2)$.  Substituting $u_{n_i}:= \left( \frac{(m-1)v_{n_i}}{m}\right)^{\frac{1}{m-1}}$ into the definition of very weak solution \eqref{def-very weak sol} and taking limits as $i\to \infty$, we also see that
\begin{equation}\label{def-U-from-V}
U(x,t) := \left( \frac{(m-1)V(x,t)}{m}\right)^{\frac{1}{m-1}}
\end{equation}
is a very weak solution of \eqref{inhomo-PME}. Equivalently, $V$ is a very weak solution of \eqref{v-PME} for all $t\in \R$, and it is classical in $Q$. This implies that $B(t)\in C^1(\R)$, it is the free boundary of $V(\cdot,t)$, and so it satisfies the Darcy's law:
\begin{equation}\label{B'=-Vx}
B'(t) = -V_x(B(t)-0,t) \geq \delta^* := -\delta'_0(0-0)>0,\qquad t\in \R.
\end{equation}

\medskip
\noindent
\underline{\it Step 3. To show that the convergence \eqref{Vni-to-V-C2} holds for the whole sequence $\{v_n\}$}. We only need to show that, if
\begin{equation}\label{vn'j-to-hatV}
v_{n'_j}(x,t)\to \widehat{V}(x,t) \mbox{ as }j\to \infty,\qquad \mbox{in the topology of } C^{2+\alpha,1+\alpha/2}_{loc}(Q),
\end{equation}
then $V(x,t)\equiv \widehat{V}(x,t)$. By contradiction, assume that, for some $s_1\in \R$,
$V(x,s_1)\not\equiv \widehat{V}(x,s_1)$. We choose a positive integer $k$ such that $-kT<s_1$.
In \eqref{exact-greater} we rewrite $x$ as $x+n$ and take $y=n-k$ for large $n$, then we have
$$
 v(x+n, t_n+(t_{n-k}-t_n);0)\geq v(x+n+1, t_{n+1}+(t_{n-k+1}-t_{n+1});0),\qquad x\in \R,
$$
that is,
$$
\tilde{v}_n(x) :=v_n(x,t_{n-k}-t_n) \geq \tilde{v}_{n+1}(x) :=v_{n+1}(x,t_{n-k+1}-t_{n+1}),\qquad x\in \R.
$$
This means that the sequence $\{\tilde{v}_n(x)\}$ is decreasing in $n$, and so it converges pointwisely as $n\to \infty$ to some function $\widetilde{V}(x)$.
On the other hand, the limits in \eqref{Vni-to-V-C2} and \eqref{vn'j-to-hatV} imply that
$$
V(x,-kT) = \widetilde{V}(x) =\widehat{V}(x,-kT),\qquad x\in \R.
$$
So, $V(x,t-kT) \equiv \widehat{V}(x,t-kT)$ for all $t>0$ since they are both solutions of
\eqref{v-PME} with the same initial data $\widetilde{V}(x)$. This leads to a contradiction when we take $t=kT +s_1>0$. This prove that the convergence \eqref{Vni-to-V-C2} holds for the whole sequence $\{v_n\}$.
Consequently, the convergence in \eqref{Bni-to-B} can also be improved to
\begin{equation}\label{Bn-to-B}
b_n(t)\to B(t) \mbox{ as }n\to \infty,\qquad \mbox{in the topology of }C_{loc}(\R).
\end{equation}

In addition, since $v(x,t;0)$ starts from $H(-x)q_2(x)$, for any $t>0$, $v(x,t;0)$ is {\it steeper than } any entire solution of \eqref{v-PME} lying in the range $[0,q_2(x)]$ (in the sense of \cite{DGM}), so does $v_n(x,t)$. In particular, $v_n(x,t_1)$ is steeper than $V(x,t_2)$. As a consequence, $V(x,t_1)$ and $V(x,t_2)$ are steeper than each other.

\medskip
\noindent
\underline{\it Step 4. To show that $V$ is a periodic traveling wave}. For any $x,t\in \R$ and any large integer $n$, we have
$$
v_{n+1}(x-1,t) =  v(x+n, t+t_n +(t_{n+1}-t_n);0) = v_n(x,t+t_{n+1}-t_n).
$$
Taking limit as $n\to \infty$ in both sides we conclude
$$
V(x-1,t)=V(x,t+T),\qquad x,\ t\in \R.
$$
This means that $V$ is a periodic sharp traveling wave of \eqref{v-PME}.

\medskip
\noindent
\underline{\it Step 5. To show that $V(x,t)\to q_0(x)$ as $x\to -\infty$ for some periodic stationary solution $q_0(x)$ of} \underline{\it \eqref{v-PME}.} We first prove a claim:
\begin{equation}\label{Vt>0}
V_t(x,t)>0,\qquad x,\ t\in \R.
\end{equation}
Differentiating $V(B(t),t)\equiv 0$ we have $V_t(B(t),t)= V_x^2(B(t)-0,t)>0$ by \eqref{B'=-Vx}. So, in order to prove the conclusion we need only to show that, for each $s_1\in \R$, $V_t (x,s_1)>0$ for all $x \in \R$.
Assume by contradiction that, for some $s_1$, there exist $x_1\in \R$ such that $V_t(x_1, s_1)\leq 0<V_t(B(s_1),s_1)$. Then, for small $s>0$, $V(\cdot, s_1+s)$ intersects $V(\cdot,s_1)$ at least once. Since they are steeper than each other (see the end of Step 3), we actually have $V(x,s_1)\equiv V(x,s_1 +s)$, and so $V_t(x,s_1)\equiv 0$, a contradiction.

Next, for any positive integer $k$, we set
$$
V_k (x,t) := V(x-k,t) \leq q_2(x),\qquad x,\ t\in \R.
$$
Since $v(x,t;1)\equiv v(x-1,t;0)\succ v(x,t;0)$, we have $v_n(x-1-k,t)\succ v_n(x-k,t)$ and so $V_{k+1}(x,t):= V(x-k-1,t)\geq V_k(x,t) := V(x-k,t)$. Thus, there exists a positive entire solution $W(x,t)$ such that
$$
V_k(x,t) \nearrow W(x,t) \mbox{ as }k\to \infty,\qquad \mbox{in the topology of } C^{2+\alpha,1+\alpha/2}_{loc} (\R^2).
$$
This also means that $V_{k+1}(x,t)\equiv V_k(x-1,t)\to W(x-1,t)$ as $k\to \infty$, and
$W(x,t)\geq \delta_4$ for all $t\in \R$ and $x\in [0,1]$. So
\begin{equation}\label{Wx-1=Wx}
W(x-1,t)\equiv W(x,t) \geq \delta_4.
\end{equation}

On the other hand, since $V$ is a periodic traveling wave, we have
$$
V_k (x-1,t) \equiv V (x-k-1,t) \equiv V(x-k, t+T) \equiv V_k(x,t+T).
$$
Taking limit as $k\to \infty$ we have $W(x-1,t)\equiv W(x,t+T)$. Combining together with \eqref{Wx-1=Wx} we have $W(x,t)\equiv W(x,t+T)$. Finally, the inequality \eqref{Vt>0} implies that  $W_t\geq 0$. So, $W_t(x,t)\equiv 0$, and $W(x,t)\equiv q_0(x)\geq \delta_4$ for some $1$-periodic stationary solution $q_0(x)$ of \eqref{v-PME}.

\medskip
\noindent
{\it Proof of Theorems \ref{thm:main1} and \ref{thm:main-limit}}. The conclusions in these two theorems  follow from the construction of $V$ and the properties we proved in the above steps.
\hfill \qed

\medskip
\begin{remark}\label{rem:ss-terrace}\rm
Note that for any give $t\in \R$ and any large $n$, $v_n(x,t)\to q_2(x)$ as $x\to -\infty$. However, the convergence $v_n\to V$ discussed in Step 3 is taking in the topology of $C^{2+\alpha,1+\alpha/2}_{loc}(\R^2)$.  Hence, the limit $q_0(x)$ of $V(x,t)$ as $x\to -\infty$  is not bigger than $q_2(x)$ but not necessarily to be $q_2(x)$. For example, in the multistable cases, $v_n$ maybe characterized by a propagating terrace, which is a combination of some traveling waves with different heights.  $V$ is actually the lowest one of them.
As we mentioned in Section 1, if we take the renormalization by considering the $h$-level set for some $h>0$, then the limit of the renormalized sequence must be a non-trivial one, and so the a priori estimates in Subsection 3.2 (which are crucial in the current paper) can be simplified in some sense. However, the renormalization limit obtained in such a way maybe a traveling wave above the second floor in the terrace, and so it is not the lowest sharp one with a free boundary.
\end{remark}

\section{$L^\infty$ Convergence in Monostable, Bistable and Combustion Equations}

We prove the $L^\infty$ convergence of $v_n$ to $V$ in monostable, bistable and combustion equations.

\medskip
\noindent
{\it Proof of Theorem \ref{thm:L-infty-cov}}. By the assumption, $p(x)$ is the minimal stationary solution of \eqref{inhomo-PME} in the range $[\kappa_0, \kappa^0]$. It is nothing but $p_1(x)$ in Proposition \ref{prop:ppss}. Denote by $q(x):=\frac{m}{m-1}p^{m-1}(x)$ the corresponding stationary solution of \eqref{v-PME}. We will prove the conclusion by showing
\begin{equation}\label{vn-to-V-Linfty}
\|v_n(x,t)-q(x)\|_{L^\infty(\R)} \to 0 \mbox{ as }n\to \infty,
\end{equation}
locally uniformly in $t\in \R$, where $v_n(x,t):= v(x+n,t+t_n;0)$ for the solution $v(x,t;0)$ of \eqref{uk-p-v} with $k=0$ and $q_2(x)\equiv q(x)$.

Under each of the conditions (i)-(iii), we see that
\eqref{F1} and \eqref{F2} hold. Then the conclusions in Theorems \ref{thm:main1} and \ref{thm:main-limit}
hold. Given $T>0$. For any small $\varepsilon>0$, by Theorem \ref{thm:main1} (ii), there exists $L_\varepsilon >0$ such that
$$
0\leq q(x)-V(x,t) \leq \varepsilon, \qquad x\in J(t): = (-\infty, B(t)-L_\varepsilon],\ t\in [-T,T].
$$
By the convergence in Theorem \ref{thm:main-limit}, we can choose $N_\varepsilon$ large such that when $n\geq N_\varepsilon$ we have
$$
\|v_n(x,t) -V(x,t)\|_{L^\infty ([B(t)-L_\varepsilon, \infty))}
 = \|v_n(x,t) -V(x,t)\|_{L^\infty ([B(t)-L_\varepsilon, B(t)+1])} \leq \varepsilon,\qquad t\in [-T,T].
$$
On the other hand, when $t\in [-T,T]$ and $n\geq N_\varepsilon$, we have by $v_n \searrow V$ pointwisely that
\begin{eqnarray*}
\|v_n(x,t) - V(x,t)\|_{L^\infty (J(t))} & \leq &
\|v_n(x,t) - q(x)\|_{L^\infty (J(t))} + \|q(x) - V(x,t)\|_{L^\infty (J(t))} \\
& \leq &  2 \|q(x) - V(x,t)\|_{L^\infty (J(t))} \leq 2 \varepsilon.
\end{eqnarray*}
This proves \eqref{vn-to-V-Linfty}, and so Theorem \ref{thm:L-infty-cov} is proved.
\hfill \qed

%
%
%

\end{document}